**Active learning-based Bayesian optimization in the realm of copper slag-blended cement systems**

Debadri Som[1,2], Rayna Maheshwari[3], Mathijs Schuurmans[4], Panagiotis Patrinos[4], Yiannis Pontikes[1]

[1]Department of Materials Engineering, KU Leuven

Kasteelpark Arenberg 44, 3001 Leuven, Belgium.

[2]RISE Research Institutes of Sweden, Department Infrastructure and Concrete Construction, Brinellgatan 4, 504 62 Borås, Sweden*

[3]Department of Materials Science and Engineering, Indian Institute of Technology (IIT) Delhi

Hauz Khas, New Delhi – 110016, India.

[4]Department of Electrical Engineering (ESAT), KU Leuven

Kasteelpark Arenberg 10, 3001 Leuven, Belgium.

**Abstract**

Accelerated mix design optimization is critical for deploying low-carbon supplementary cementitious materials (SCMs) because traditional experimental approaches require extensive testing campaigns. The authors demonstrate that Bayesian optimization (BO) can identify near-optimal blended cement formulations using an AI-driven approach with minimal experimental data. Starting with only 10 initial experiments and a 2:1 data-to-variable ratio reflecting realistic laboratory constraints, the authors optimized copper slag-limestone-Portland cement systems for 32.5N strength class with respect to 2-day compressive strength, cost, and $CO_2$ emissions. Gaussian process surrogate models guided sequential experimentation through Expected Improvement and Upper Confidence Bound acquisition functions. Within 2-6 iterations, BO identified Pareto non-dominated solutions meeting strength requirements (≥10 MPa at 2 days, ≥32.5 MPa at 28 days) while achieving $CO_2$ emissions below 500 kg $CO_2$/ton; consistent with industry decarbonization targets. A conservative update strategy incorporating uncertainty bounds for 28-day strength enabled rapid iteration without waiting 28 days per cycle. Comparative analysis revealed that linear kernels outperformed nonlinear alternatives in predictive accuracy, though radial basis function kernels were preferred for active learning due to superior uncertainty quantification. This work demonstrates BO as a practical decision-support tool for cement research under severe data constraints.

**Keywords**

Bayesian optimization, AI-driven material design, blended cements, Fe-rich slags, multi-objective optimization

**Abbreviations**

| | |
|---|---|
| BO | Bayesian Optimization |
| $\hat{c}(x)$ | Deterministic cost of cement mix design (x), computed from material proportions and unit costs |
| $c_{\max}$ | Maximum allowable cost constraint |
| $\mathcal{C}$ | Set of box and linear constraints defining feasible designs |
| $\mathcal{D}_0$ | Initial dataset of experimentally evaluated designs |
| $\mathcal{D}_r$ | Ranked feasible designs |
| $\mathcal{D}_t$ | Dataset available at Bayesian optimization iteration (t) |
| $\hat{e}(x)$ | Deterministic $CO_2$ emissions of cement mix design (x), computed from material proportions and emission factors |
| $e_{\max}$ | Maximum allowable $CO_2$ emissions constraint |
| $\exp(\cdot)$ | Exponential function |
| $\mathbb{E}[\cdot]$ | Expectation operator |
| EI | Expected Improvement acquisition function |
| $\mathcal{F}$ | Set of feasible candidate designs after applying constraints |
| $f_2(x)$ | Gaussian Process function modelling 2-day compressive strength |
| $f_{28}(x)$ | Gaussian Process function modelling 28-day compressive strength |
| $\Gamma(\cdot)$ | Gamma function |
| $\mathcal{GP}$ | Gaussian Process |
| $\kappa$ | Tunable parameter for UCB acquisition function (kappa) |
| $k_2(x, x')$ | Kernel (covariance) function for the 2-day Gaussian Process |
| $k_{28}(x, x')$ | Kernel (covariance) function for the 28-day Gaussian Process |
| $K_\nu(\cdot)$ | Modified Bessel function of the second order of $\nu$ |
| $\mu(x)$ | Gaussian Process posterior mean |
| $\mu_2(x)$ | Gaussian Process posterior mean prediction of 2-day strength at (x) |
| $\mu_{28}(x)$ | Gaussian Process posterior mean prediction of 28-day strength at (x) |
| $\tilde{\mu}_2(x)$ | Normalized 2-day compressive strength (mean from Gaussian process) |
| $\mu_{2,\max}$ | Maximum 2-day compressive strength in the dataset |
| $\mu_{2,\min}$ | Minimum 2-day compressive strength in the dataset |
| $\nu$ | Smoothness parameter of the Matérn kernel |
| $\mathcal{N}$ | Number of candidate designs sampled per Bayesian optimization iteration |
| $\Phi(.)$ | Standard cumulative distribution function (capital phi) |
| $\phi(.)$ | Standard probability density function (lowercase phi) |
| $r(a, a')$ | Scaled Euclidean distance between two input points in a kernel |
| $R^2$ | Coefficient of determination |
| $s(x)$ | Scalarized objective value at design (x): $(s(x) = w_1\mu_2(x) - w_2\hat{c}(x) - w_3\hat{e}(x))$ |
| $s^*$ | Best scalarized objective value observed so far |
| SD | Standard Deviation |
| $\sigma_2(x)$ | Gaussian Process posterior standard deviation of 2-day strength at (x) |
| $\sigma_{2,\min}$ | Minimum required compressive strength at 2 days (MPa) |
| $\sigma_{28}(x)$ | Gaussian Process posterior standard deviation of 28-day strength at (x) |
| $\sigma_{28,\min}$ | Minimum required compressive strength at 28 days (MPa) |
| $\tilde{\sigma}_s(x)$ | Scalarized standard deviation used in acquisition function |
| $t$ | Bayesian optimization iteration index |
| $T$ | Final iteration index at termination |
| UCB | Upper Confidence Bound |
| $w = (w_1, w_2, w_3)$ | Weights used in scalarized objective function |
| $x$ | Input vector representing a cement mix design |
| $x'$ | Another input vector in the design space (used in kernel definition) |
| $x^*$ | Selected candidate design maximizing the acquisition function |
| $\mathcal{X}$ | Design space of admissible cement mix designs |
| $y$ | Vector of 2-day compressive strength, cost and $CO_2$ equivalent |
| $y_2$ | Observed 2-day compressive strength |
| $y_{28}$ | Observed 28-day compressive strength |

$z$ z-score, defined as $\frac{s(x)-s^*}{\sigma_s(x)}$

## 1. Introduction

Portland cement accounts for approximately 8% of anthropogenic $CO_2$ emissions (Ellis *et al.*, 2020). This has propelled recent research directions to the area of novel Supplementary Cementitious Materials (SCMs)(Snellings, Suraneni and Skibsted, 2023), especially when the availability of common SCMs like fly ash is expected to reduce (International Energy Agency, 2009), owing to the prevalence of alternative energy sources. Non-ferrous metallurgical slag (NFMS) are typically Fe-rich residues from primary or secondary production of metals like copper. In recent studies, they have shown potential as a prospective SCM, albeit with differences in chemical composition from typical SCMs like fly ash and ground granulated blast furnace slag (GGBFS). Hallet *et al.* (2022) used a laboratory-treated Fe-rich slag with successful replacements up to 30 wt%, reaching a strength class of 52.5N as per EN 197-1 (European Committee for Standardisation, 2011). However, such slags have shown to reduce early-age strength which is a key limitation of Fe-rich slags (Hallet, De Belie and Pontikes, 2020). Different hypotheses have been proposed for the negative impact of Fe-rich slag on early-age strength including the reduction of nucleation sites for calcium silicate-hydrate (C-S-H), leaching of elements like Zn (Hallet, De Belie and Pontikes, 2020) and the presence of $Fe^{2+}$ which affects the hydration of alite and belite phases in a blended cement system. Similar findings corresponding to lower early-age strength of Fe-rich slag (copper slag)-based binder systems were reported by Feng *et al.* (2019) and Mahajan and Muhammad (2024). A key strategy to negate the effect of reduced early-age strength is the addition of limestone powder, creating a ternary Fe-rich slag-limestone-Portland cement blended cement system (Hallet *et al.*, 2023; Giels *et al.*, 2025). This is typically because limestone provides more nucleation sites for C-S-H, contributing to early-age strength. Additionally, limestone facilitates the formation of phases like hemicarboaluminate and monocarboaluminate, preventing the formation of monosulphoaluminates from ettringite, which is critical for early age compressive strength. Hence, the use of ternary systems with Fe-rich slag, limestone and Portland cement has shown to be an effective $CO_2$ equivalent reduction strategy. However, despite promises, such 'novel' SCMs are extremely slow to penetrate the market with a typical adoption time for new cementitious materials around 10 – 20 years (US Department of Energy, 2023). This is generally due to extended periods of testing of cementitious systems which eventually delays the publication of standards. The United States Department of Energy recommends investments in accelerated testing of materials for quicker market adoption of novel SCMs (US Department of Energy, 2023). The use of effective data-driven techniques can be a game-changing strategy in this regard with users being able to take smarter and quicker decisions based on trends suggested by historical data.

Data-driven studies in the field of building materials have generally been dominated by predictive machine learning algorithms like linear regression (Tam, 2022), random forest regression (Chen *et al.*, 2025; Mesfin and Kim, 2026), support vector machines (Deng *et al.*, 2024) and artificial neural networks (Piro, Mohammed and Hamad, 2023). While these are all powerful tools to get predictive insights from historical data, most of them require a high data-to-independent variable ratio which is often challenging to obtain in the field of building materials (Li *et al.*, 2022). Furthermore, concerning data-driven optimization of mortar and concrete mix designs, the emphasis has typically been on metaphor-inspired metaheuristic algorithms (Lee, Yoon and Kim, 2012; Huang *et al.*, 2020; Afshan and Ramezanianpour, 2026). However, as argued by Sörensen (2015), majority of these metaheuristic algorithms are redundant in terms of contributing to the scientific literature for solving non-convex problems. From a mathematical standpoint, convex optimization methods such as linear programming and quadratic programming, are preferable whenever the problem structure permits, as they provide strong theoretical guarantees of global optimality (Boyd and Vandenberghe, 2004). Nevertheless, many practically relevant problems in the optimization of building materials' properties are non-convex which do not offer straightforward global optimality guarantees compared to their convex counterpart. Furthermore, extensions to uncertainty-aware techniques using methods like distributionally robust optimization (Schuurmans and Patrinos, 2023) enable the consideration of robustness while solving optimization problems. Historically, the domain of optimization in building materials has relied on single objectives which are inherently easier to solve (Roy *et al.*, 2024). However, conflicting objectives (for example, compressive strength and cost) can warrant the need for efficient optimization techniques like desirability method (Heidari, Rivard and Wilson, 2024), scalarization (Bharadwaj, Isgor and Weiss, 2024) and Pareto trade-offs (Qiu *et al.*, 2025; Mumtaz *et al.*, 2026; Xu, Chen and Liu, 2026). While such techniques provide a basis for optimizing difficult multi-objective situations, they do not sequentially learn from historical data. Sequential learning can significantly reduce the number of required data points which is critical in a slow-progressing field like that of building materials.

Formally, Bayesian optimization (BO) is one of the most popular cases of active learning, which further, is a subset of sequential learning. The term 'Bayesian optimization' originates from a Bayesian-like approach to optimization based on posterior-prior relationships (Brochu, Cora and Freitas, 2010). Over the last two decades, BO has been rigorously applied to a wide range of studies in materials science like alloy discovery based on CALPHAD simulations (Sundar *et al.*, 2025), studying properties of polyamide nanocomposites (Ozdemir *et al.*, 2025), development of microstructure-based structure-property models (Liu *et al.*, 2024) and others. Kalidindi (2019) attributes the rise of using BO techniques in materials science to the Material Genome Initiative (MGI) (Materials Genome Initiative (U.S.), 2011) which has been instrumental in cross-disciplinary progress of scientific research especially bridging data-driven research and materials science. On a similar note, Lookman *et al.* (2019) highlighted the potential of using BO techniques for material discovery especially with the prevalence of high-throughput experimental techniques. Further recent work on BO has focussed on the development of tools like BoTorch which is a modular tool for the application of BO to real-world problems (Balandat *et al.*, 2020). However, on the front of building materials, the application of BO for the discovery of new mix formulations has been rather scarce. Saleh *et al.* (2022) were one of the first to use the concept of adaptive design in the context of building materials, specifically ultra-high performance concrete (UHPC). However, they only focused on a single-objective Bayesian optimization with 28-day compressive strength as the single target. In order to have a holistic understanding of the performance, economics and environmental dimensions of binder/mortar/concrete systems, there have been multiple studies on multi-objective BO in the field of building materials (Ament *et al.*, 2023; Völker *et al.*, 2023; Pfeiffer *et al.*, 2024). However, the works of Ament *et al.* (2023) and Pfeiffer *et al.* (2024) used more than thousand data points as a starting point for BO which can be impractical for adaptive experimental design in a laboratory setting. Völker *et al.* (2023), on the other hand, randomly sampled between four and 20 data points from a database to study the performance of Gaussian process models with less training data, which is a more practical application of BO in the field of cement and concrete research.

Since most BO studies employ a small initial dataset, the design of the initial dataset is critical in this context. Völker *et al.* (2023) employed a related strategy for training GP regression models in a Bayesian optimization context, where multiple combinations of four to 20 data points were evaluated. Their results demonstrated a clear improvement in GP predictive performance as the number of training points increased. However, in that study the training datasets were sampled randomly, without explicit consideration of balanced or near-orthogonal coverage of the input space. As a result, comparisons across different sampled datasets may be influenced by variations in input-space coverage, complicating direct performance comparisons. On a similar note, Jones and Schonlau (1998) recommended space filling designs as a starting point for computer experiments. While that is a recommended option for cases with easy-to-collect data/simulation results, it is not feasible to use such techniques in the field of building materials as a starting data set for sequential experiments.

Given the recent advances of using active learning techniques in the field of building materials, it is necessary to work on best practices in this domain to accelerate the development of low-carbon mix designs. This paper aims to use a small starting dataset on copper slag-based ternary blended cement systems and develop a methodology to predict the next experiments. Furthermore, this paper generates simulated mix designs and their corresponding compressive strength to augment the initial dataset and converge towards the global Pareto front. Unlike previous studies, this work focuses on BO as a practical day-to-day experimental decision-making tool under realistic laboratory constraints, rather than as a purely retrospective data-mining approach.

In this paper, the authors address several key questions that extend existing literature at the interface of applied mathematics and building materials. First, they investigate whether an orthogonal or near-orthogonal design can serve as an effective starting dataset for active learning, particularly under constrained conditions where the data-to-variable ratio is as low as 2:1. This is motivated by prior studies, such as Nahvi *et al.* (2019), who employed a factorial design with 16 data points for two variables with four levels, and Völker *et al.* (2023), who observed a reduction in optimization efficiency when the dataset size was decreased from 20 to 4, potentially due to increased noise in literature-derived data. In this context, the influence of categorical features on model performance is also examined. Second, the study explores how kernel selection, length-scale definition, and predictive uncertainty in compressive strength models can be systematically incorporated into the Bayesian optimization (BO) framework. Third, recognizing that 28-day compressive strength is the most relevant performance metric for mortar mix design but introduces significant delays in iterative experimentation, the authors propose strategies to address this limitation within the optimization loop. Finally, from an application perspective, the work demonstrates how a hybrid design approach, combining an initial near-orthogonal design with sequential active learning can be used to generate blended cement mix designs that simultaneously satisfy target strength classes, $CO_2$ equivalent constraints defined by CEMBUREAU, and user-specified cost limits. This work adds on to the previous work of the authors (Som *et al.* (2026)) where it provides a small-data alternative and the possibility to do sequential experiments rather than the traditional design-experiment-analyse approach, reaching optimal mix designs with fewer experimental trials.

## 2. Methods

### *2.1 Description of the dataset*

The dataset used in this study is based on an Fe-rich slag-based blended cement system with limestone and Portland cement. The Fe-rich slag is an industrial slag from Aurubis, Hamburg while limestone and Portland cement were sourced from Carmeuse, Belgium and Holcim, respectively. The original dataset consists of 51 data points which were synthesized using a 3-stage experimental design combining definitive screening design (DSD) (Jones and Nachtsheim, 2011) and fast flexible filling design (FFFD) which is a space-filling design. Som *et al.* (2026) shows the full dataset related to the problem. Table 1 shows a summary of the different independent variables used in the study and their levels/ranges. The dependent variables were 2-day and 28-day compressive strength, cost and $CO_2$ equivalent. Further details regarding the modelling of the response variables have been shown in section 2.2.

**Table 1**: Details of the independent variables in this paper

| Factor | Type of variable | Level 1 | Level 2 |
|---|---|---|---|
| Copper slag | Continuous | 20 wt % of total binder | 40 wt % of total binder |
| Limestone-cement ratio | Continuous | 0.15 | 0.55 |
| Water-binder (wb) ratio | Continuous | 0.45 | 0.50 |
| Cement type | Discrete/Categorical | CEM I 52.5 N | CEM I 52.5 R |
| Fe-rich slag fineness | Discrete/Categorical | 2700 $cm^2/g$ | 4400 $cm^2/g$ |

From the full dataset of 51 data points, a subset of 10 data points was used in this study. This choice was made to highlight the ability of Bayesian optimization to operate effectively in small-data regimes and to propose informative subsequent experiments. Such a setting closely reflects typical laboratory conditions in cement and concrete research, where developing and testing many mix designs, particularly at 28 days, is often costly in terms of time, materials, and manpower. Consequently, experimental studies frequently rely on limited datasets to guide further experimentation toward optimal material performance.

The 10 data points used in this study are listed in Table 2. These points were selected to ensure maximal diversity across both continuous and categorical variables, with the aim of providing broad coverage of the design space. In principle, a small-scale orthogonal or fractional factorial design (e.g., a quarter-fractional factorial design) would be advantageous, as it ensures balanced representation of factor levels for the Gaussian process (GPs) models to learn without bias. The used dataset is the closest to such a situation that could be obtained from the full dataset (Som *et al.*, 2026). Furthermore, the categorical variables in the feature set were one-hot encoded before developing the surrogated models. To estimate the kernel performance, a validation dataset was also used which was randomly sampled (with fixed seed) from the 51-point dataset. However, the training data points were avoided in that case to ensure the model has not encountered those points before while training.

**Table 2**: Training dataset for the GP models used in this paper

| Fe-rich slag | Limestone | Cement | wb ratio | Copper slag fineness | Cement type | 2-day strength | 28-day strength |
|---|---|---|---|---|---|---|---|
| 20 | 10.4 | 69.6 | 0.45 | 2700 | R | 29.6 | 52.7 |
| 20 | 10.4 | 69.6 | 0.5 | 2700 | N | 21.8 | 48.1 |
| 20 | 28.4 | 51.6 | 0.45 | 4400 | N | 16.3 | 39.4 |
| 20 | 28.4 | 51.6 | 0.45 | 4400 | R | 15.4 | 39.0 |
| 40 | 7.8 | 52.2 | 0.45 | 4400 | N | 14.6 | 43.0 |
| 40 | 7.8 | 52.2 | 0.5 | 2700 | N | 14.6 | 27.6 |
| 40 | 21.3 | 38.7 | 0.45 | 2700 | N | 9.3 | 25.4 |
| 40 | 21.3 | 38.7 | 0.5 | 2700 | R | 9.8 | 25.1 |
| 20 | 28.4 | 51.6 | 0.5 | 4400 | R | 13.9 | 33.4 |
| 20 | 10.4 | 69.6 | 0.5 | 2700 | R | 25.6 | 52.6 |

### *2.2 Case study*

In this paper, the technique of BO-based active learning has been demonstrated on a specific case of 32.5N cement strength class (European Committee for Standardisation, 2011). CEMBUREAU (the European cement association) has formally committed to decarbonization aligned with the Paris Agreement, framed in its 2050 Net Zero Roadmap (CEMBUREAU, 2020), which emphasizes the introduction and market deployment of low-carbon cement and concrete solutions to reduce lifecycle $CO_2$ emissions. This strategic focus inherently supports the development and wider use of cement classes that allow high replacement of clinker with supplementary materials, such as the 32.5 N strength class, over higher strength classes that typically require higher clinker content (especially in the case of low-reactive SCMs like Fe-rich slags). Furthermore, 32.5N cement strength class relaxes constraints on compressive strength (as compared to 42.5N and 52.5N cement strength class) which increases the possibility of a higher degree of Portland cement replacement with limestone and Fe-rich slag. The relaxation of the constraint on compressive strength enables a clearer demonstration of BO in terms of cement clinker replacement.

#### *2.2.1 Optimization problem formulation*

The cost and the equivalent $CO_2$ models were assumed to be linear models, as a function of the unit price of the raw materials and the weight percentage of their presence in the binder system out of 100 wt% binder. Equations 1 and 2 show the mathematical expression for the total cost and total $CO_2$ of the blended cement systems while further information on unit cost and unit equivalent $CO_2$ can be found it Som *et al.* (2026).

$$\hat{c}(x) = 0.01(c_{ce}c_c(x) + c_{ls}c_l(x) + c_{sl}c_s(x)) \quad (1)$$

$$\hat{e}(x) = 0.01(e_{ce}c_c(x) + e_{ls}c_l(x) + e_{sl}c_s(x)) \quad (2)$$

where, $\hat{c}(x)$ and $\hat{e}(x)$ are total cost and total equivalent $CO_2$ calculated for a set of independent variables $x$ while $c_{ce}, c_{ls}$ and $c_{sl}$ are the unit costs of Portland cement, limestone and Fe-rich slag respectively and $c_c, c_l$ and $c_s$ are the percentages of Portland cement, limestone and Fe-rich slag in the total binder system (constrained to 100 wt%). In Equation 2, $e_{ce}, e_{ls}$ and $e_{sl}$ are the unit equivalent $CO_2$ values of Portland cement, limestone and Fe-rich slag respectively.

#### *2.2.2 Objective function*

The objective function for the multi-objective problem in this paper was built based on scalarization. However, since compressive strength [MPa], cost [EUR/ton] and equivalent $CO_2$ [kg $CO_2$/ton] have different units and different ranges of values, min-max normalization was used to have them in the range of 0 to 1. The normalization equation is shown in Equation 3.

$$Norm\ (y) = \frac{y - \min(y)}{\max(y) - \min(y)} \quad (3)$$

where, $Norm(y)$ is the normalized form of scalar objective value y, $min(y)$ and $max\ (y)$ are the global minimum and global maximum values obtained from the initial training dataset.

In the case of this paper, 2-day compressive strength has been modelled using GPs while cost and equivalent $CO_2$ models have been based on Equations 1 and 2. This means that that the scalarized model is also a GP with a different mean but the same variance. Equations 4 - 6 show the scalarized mean and standard deviation.

$$s(x) = w_1\tilde{\mu}_2(x) - w_2\tilde{c}(x) - w_3\tilde{e}(x) \quad (4)$$

$$\tilde{\sigma}_s(x) = \frac{w_1\sigma_2(x)}{\mu_{2,max} - \mu_{2,min}} \quad (5)$$

$$w_1 + w_2 + w_3 = 1 \quad (6)$$

In Equations 4 and 5, $s(x)$ represents the scalarized mean objective derived from the Gaussian process (GP) models. The weights $w_1$, $w_2$, and $w_3$ correspond to the relative contributions of 2-day compressive strength, cost, and equivalent $CO_2$ emissions, respectively. The quantities $\tilde{\mu}_2(x), \tilde{c}(x)$ and $\tilde{e}(x)$ denote the predicted normalized mean responses obtained from the 2-day compressive strength GP model, the deterministic cost model, and the deterministic equivalent $CO_2$ model, respectively, while $\tilde{\sigma}_s(x)$ denotes the predictive standard deviation of the 2-day compressive strength GP model.

For the scalarized mean objective, all three responses contribute through their respective weights, which were set equal to 0.33 ($w_1 = w_2 = w_3$) in the present study, thereby assigning equal importance to early-age strength, cost, and equivalent $CO_2$ emissions. In contrast, the scalarized standard deviation is solely determined by the uncertainty associated with the 2-day compressive strength prediction. This is because both the cost and equivalent $CO_2$ models are deterministic and therefore do not contribute to predictive uncertainty in the scalarized objective formulation.

#### *2.2.3 Constraints*

The use of constraints in an optimization problem restricts solutions to practical feasibility and the original data structure. Based on the problem in question, different classes of constraints were considered in this problem which have been discussed in this subsection.

- *Material constraints/Box constraints*

Equation 7 represents the summation constraint for the total wt% of binders. This must be followed to keep the total binder content constant which affects the density, flow and other properties of the mortar mix. The notations $c_c$ , $c_l$ and $c_s$ are the wt% of cement, limestone and Fe-rich slag in a binder system.

$$c_c + c_l + c_s = 100 \tag{7}$$

Equations 8 - 10 represent the individual material bounds on Fe-rich slag, limestone and Portland cement. It is necessary to adhere to this constraint to prevent accidental extrapolation in the model outcomes.

$$c_{cmin} \leq c_c \leq c_{cmax} \tag{8}$$

$$c_{smin} \leq c_s \leq c_{smax} \tag{9}$$

$$c_{lmin} \leq c_l \leq c_{lmax} \tag{10}$$

where $c_{cmin}$ , $c_{smin}$ and $c_{lmin}$ are the minimum content of cement, Fe-rich slag and limestone in the binder system (corresponding to the data set) respectively and $c_{cmax}$ , $c_{smax}$ and $c_{lmax}$ are the maximum content of cement, Fe-rich slag and limestone in the binder system (corresponding to the data set) respectively.

Equation 11 shows the limit on wb ratio which was added as a constraint to ensure the flow diameter of the mortar mixes to be in the range 180 – 200 mm, for a consistent comparison of compressive strength. The total binder content in this case includes Portland cement, limestone and Fe-rich slag.

$$0.48 \leq \frac{water}{binder} \leq 0.50 \tag{11}$$

A final compositional constraint was added to adhere to the limits of limestone/cement ratio which was an initial component of the experimental design. Equation 12 shows the inequality defining the constraint imposed on limestone/cement ratio. This constraint is corresponding to the limestone content defined in CEM II/A-L and CEM II/B-L as defined in EN 197-1 (European Committee for Standardisation, 2011).

$$0.15 \leq \frac{c_l}{c_c} \leq 0.55 \tag{12}$$

- *Compressive strength constraints*

While the material constraints make sure the model results adhere to compositional constraints of the materials, compressive strength constraints are necessary for adherence to the 32.5N strength class cement system. A 32.5N cement strength class must exhibit a mortar compressive strength of at least 10 MPa at 2 days and at least 32.5 MPa at 28 days (European Committee for Standardisation, 2011) after casting. Equations 13 - 14 show the mathematical expressions for the compressive strength constraints used in this paper.

$$\mu_2(x) - 1.96\sigma_2(x) \geq 10 \text{ MPa} \tag{13}$$

$$\mu_{28}(x) - 1.96\sigma_{28}(x) \geq 32.5 \text{ MPa} \tag{14}$$

From Equations 13 and 14, the compressive strength constraints were formulated using lower confidence bounds of the GP predictions. Specifically, the terms 1.96 $\sigma_2(x)$ and 1.96 $\sigma_{28}(x)$ were subtracted from the corresponding predictive means, such that the lower bounds of the 95% confidence intervals were required to exceed 10 MPa and 32.5 MPa at 2 and 28 days after casting, respectively. This formulation penalizes candidate mix designs associated with high predictive uncertainty and ensures that the strength requirements are satisfied with high statistical confidence. From an engineering perspective, the terms 1.96 $\sigma_2(x)$ and 1.96 $\sigma_{28}(x)$ can be interpreted as safety margins that reduce the risk of selecting uncertain designs whose true compressive strength may fall below the specified limits. Also, it is to be noted that 28-day compressive strength is only used as a constraint and is not a part of the objective function. This is because model updates for 28-day compressive strength would require waiting for 28 days, delaying the iteration loop. However, the authors do acknowledge the limitation that using 2-day compressive strength as a proxy for 28-day compressive strength in the objective function could lead to missing out promising mixes especially with an SCM in the mix design.

- *Cost and equivalent $CO_2$ constraints*

The current problem under consideration is multi-objective (compressive strength, cost and equivalent $CO_2$) and hence constraints for cost and equivalent $CO_2$ must also be added pertaining to a 32.5N cement strength class. In this case a $c_{max}$ (maximum permissible cost) value of 110 EUR/ton has been used which is highly dependent on the sourcing of raw materials. In the context of equivalent $CO_2$, an $e_{max}$ (maximum permissible equivalent $CO_2$) value of 500 kg $CO_2$/ton has been used which lesser than target equivalent $CO_2$ for 2050, reported by CEMBUREAU (CEMBUREAU, 2020), albeit with a lower strength class. Mathematically, the expressions for the cost and equivalent $CO_2$ constraints have been shown in Equations 15 - 16 where $\hat{c}(x)$ is the computed cost of the mix design (considering only binders) and $e(x)$ is the computed equivalent $CO_2$ of the mix design (considering only binders).

$$\hat{c}(x) \leq c_{\max} \tag{15}$$

$$\hat{e}(x) \leq e_{\max} \tag{16}$$

### *2.3 BO technique*

The BO technique used for sequential learning in this paper is shown in Algorithm 1. Each component of the algorithm is further discussed in this section (except points 1 and 2 which are mostly related to the treatment of data). Please refer to the set of abbreviations at the beginning of the paper on information about the symbols used.

**Algorithm 1:** Multi-objective BO through scalarization technique

**Input**: Initial dataset $\mathcal{D}_0 = \{(x_i, y_{2,i}, y_{28,i})\}_{i=1}^{n}$ ; design space $\mathcal{X}$; constraints $\mathcal{C}$; weights $w = (w_1, w_2, w_3)$; candidate size $N$ (section 2.1)

**Output**: Final dataset $\mathcal{D}_T$ and recommended designs $x*$

1. **Initialization of data set**: Set $t = 0, \mathcal{D}_t = \mathcal{D}_0$.

2. **GP models**: Define GPs for 2-day and 28-day strength with $\mu$ as the mean and $k$ as the kernel function (section 2.3.1):

$$f_2(x) \sim \mathcal{GP}(\mu_2, k_2)$$
$$f_{28}(x) \sim \mathcal{GP}(\mu_{28}, k_{28})$$

3. **Bayesian Optimization (BO) loop**: while stopping criterion not met (defined in section 2.2.3)

3.1 **Training**: Fit GPs on $\mathcal{D}_t$ to obtain $(\mu_2, \sigma_2), (\mu_{28}, \sigma_{28})$.

3.2 **Normalization**: $\tilde{\mu}_2(x) = \frac{\mu_2(x) - \mu_{2,\min}}{\mu_{2,\max} - \mu_{2,\min}}, \tilde{c}(x) = \frac{\hat{c}(x) - c_{\min}}{c_{\max} - c_{\min}}, \tilde{e}(x) = \frac{\hat{e}(x) - e_{\min}}{e_{\max} - e_{\min}}$

**3. Scalarization** (section 2.2.2): $s(x) = w_1\tilde{\mu}_2(x) - w_2\tilde{c}(x) - w_3\tilde{e}(x)$

$$\tilde{\sigma}_s(x) = \frac{w_1\sigma_2(x)}{\mu_{2,\max} - \mu_{2,\min}}$$

$$s^* = \max_{x_i \in \mathcal{D}_t} s(x_i); \text{ where } s^* \text{ is best scalar}$$

3.4 **Sampling**: Generate *N* candidates $x_j \sim \mathcal{U}(\mathcal{X})$ satisfying mass balance.

3.5 **Prediction**: Evaluate $\mu_2(x_j), \sigma_2(x_j), \mu_{28}(x_j), \sigma_{28}(x_j)$

3.6 **Acquisition** (section 2.3.2): compute Expected Improvement (EI) or Upper Confidence Bound (UCB): shown only for EI due to brevity.

$$\mathrm{EI}(x) = \left(s(x_j) - s^*\right)\Phi(z) + \sigma_s(x_j)\phi(z),\ z = \frac{s(x_j) - s^*}{\sigma_s(x_j)}, \text{ using}$$

$$s(x_j), \sigma_s(x_j)$$

3.7 **Retain $\boldsymbol{x_j, j = 1, \dots, N}$ s.t.** (section 2.2.3)

$$\mu_2(x_j) - 1.96\sigma_2(x_j) \geq \sigma_{2,min};$$

$$\mu_{28}(x_j) - 1.96\sigma_{28}(x_j) \geq \sigma_{28,min;};\hat{c}(x_j) \leq c_{\max}, \hat{e}(x_j) \leq e_{\max}, x_j \in \mathcal{C}$$

$$\mathcal{F} = \{x_j | constraints\ on\ x_j\ (section\ 2.2.1.2), j = 1, \dots, N\}$$

3.8 **Selection**: $x^* = \underset{x_j \in \mathcal{F}}{\mathrm{argmax}}\ \mathrm{EI}(x_j)$

3.9 **Evaluation**: Conduct experiment or use predicted values to obtain $x^*$ $y_2^*, y_{28}^*$

3.10 **Update** (section 2.3.4): $\mathcal{D}_{t+1} = \mathcal{D}_t \cup \{(x^*, y_2^*, y_{28}^*)\}$

$$t \leftarrow t + 1$$

4. **Return**: $\mathcal{D}_t$ and ranked feasible designs $\mathcal{D}_r$

---

### *2.3.1 Gaussian process modelling*

Referring to point 3 in the algorithm, GPs were used to train the models for 2-day and 28-day compressive strength. Gaussian Processes (GPs) are non-parametric models widely used as surrogate models for expensive-to-evaluate functions, particularly in the context of Bayesian optimization and design of experiments. A GP defines a probability distribution over functions and provides a flexible framework for modelling unknown relationships between input variables and observed responses while explicitly quantifying predictive uncertainty.

Formally, a Gaussian Process is defined as a collection of random variables, any finite subset of which follows a joint multivariate Gaussian distribution. A GP prior over a function $f(x)$ is specified as shown in Equation 17.

$$f(x) \sim \mathcal{GP}(\mu(x), k(x, x')) \tag{17}$$

where $\mu(x)$ is the mean function and $k(x, x')$ is the covariance/kernel function. $x$ and $x'$ are two input points in the design space. In most practical applications, the mean function is assumed to be zero (in this paper) without loss of generality, and prior assumptions about smoothness, variability, and correlation structure are encoded through the choice of kernel. Hence, it can be understood that the choice of kernel and its corresponding hyperparameters are critical for data fitting in the context of GP surrogates for the BO loop. The kernels used in this paper and their corresponding equations are shown in Table 3.

**Table 3**: Kernels used in this paper

| **Kernel name** | **Kernel equation** |
|---|---|
| Radial basis function (RBF) | $k(a, a') = \sigma_f^2 \exp\left(-\frac{1}{2} r^2(a, a')\right);\ r(a, a') = \sqrt{\sum_{j=1}^{d} \frac{(a_j - a_j')^2}{\ell_j^2}}$ |
| Matérn 3/2 | $k_\nu(a, a') = \sigma_f^2 \frac{2^{1-\nu}}{\Gamma(\nu)} \left(\sqrt{2\nu} r(a, a')\right)^\nu K_\nu\left(\sqrt{2\nu} r(a, a')\right); \nu = 3/2$ |
| Matérn 5/2 | $k_\nu(a, a') = \sigma_f^2 \frac{2^{1-\nu}}{\Gamma(\nu)} \left(\sqrt{2\nu} r(a, a')\right)^\nu K_\nu\left(\sqrt{2\nu} r(a, a')\right); \nu = 5/2$ |
| Linear (benchmark) | $k(a, a') = \sigma_f^2 a^\top a' + \sigma_0^2$ |

In Table 3, $a$ and $a'$denote two distinct input data points used in the kernel definition, while $k(a, a')$ represents the kernel (covariance) function that quantifies similarity between their corresponding responses. The parameter $\ell$ (or $\ell_j$ in the case of feature-dependent length scales) defines the characteristic length scale of the kernel and governs the smoothness of the Gaussian process response surface. Smaller values of $\ell$ allow rapid local variations in the predicted response and may lead to overfitting, whereas larger values enforce smoother response surfaces and may result in underfitting. From a materials perspective, the length scale can therefore be interpreted as a measure of the sensitivity of compressive strength to variations in mix composition/cement type for example. For descriptions related to the other symbols, please refer to the list of abbreviations at the beginning of the paper.

In practice, kernel length scales may be selected heuristically based on the physical ranges of the mix-design variables or inferred through maximization of the log marginal likelihood (LML). However, as discussed by Rasmussen and Williams (2006), LML-based hyperparameter optimization can be unreliable in data-limited settings due to the possible presence of multiple local minima with conflicting hyperparameters settings. Given the limited number of experimental data points available in this study, kernel hyperparameters were therefore fixed rather than optimized. Specifically, a grid-search was conducted on selected values of the length-scales (close to the range of the variable in the dataset) for the different independent variables. The selected values were as follows:

- Length-scale candidates for copper slag and limestone (wt%): 10, 20, 30
- Length-scale candidates for water-to-binder ratio: 0.01, 0.05, 0.1, 0.5
- Length-scale candidates for cement type and Fe-rich slag fineness: 1, 2, 3

The performance of the kernels was judged based on test $R^2$ on a held-out validation set consisting of 10 data points extracted from the original 51-point dataset using a fixed random seed to ensure repeatability.

While this approach does not guarantee globally optimal hyperparameters, it provides a stable and conservative modelling strategy that is well suited to experimental cement and concrete research, where data acquisition is costly and sample sizes are necessarily small. The adopted strategy prioritizes robust uncertainty estimation over aggressive fitting, which is particularly important for guiding sequential experimental design.

As shown in Table 3, $\sigma_f^2$ denotes the variance of the Gaussian process function and was fixed to unity in this work, consistent with the normalization of the response variables. In the case of the linear kernel, $\sigma_0^2$ represents a bias (offset) term. The linear kernel was included as a benchmark to provide a baseline level of predictive performance, against which the nonlinear kernels were evaluated. The remaining kernels considered: Matérn 3/2, Matérn 5/2, and the radial basis function (RBF) kernel are nonlinear, anisotropic kernels with controllable smoothness and are among the most commonly employed kernel/covariance functions in Gaussian process regression and Bayesian optimization (Rasmussen and Williams, 2006). Notably, the RBF kernel can be interpreted as a limiting case of the Matérn family corresponding to infinitely smooth response surfaces.

### *2.3.2 Acquisition function selection*

In BO, acquisition functions play a central role by guiding the sequential selection of new experiments based on the predictions of the surrogate model. While the Gaussian process provides a probabilistic estimate of the response surface in terms of a predictive mean and uncertainty, the acquisition function translates this information into a quantitative criterion that balances exploration of uncertain regions of the design space with exploitation of regions predicted to yield high performance. At each iteration, the acquisition function is optimized to identify the most informative next experiment, thereby enabling efficient learning from a limited number of costly evaluations.

In this paper, Expected Improvement (EI) and Upper Confidence Bound (UCB) are the two acquisition functions employed, both of which are among the most commonly used strategies in the Bayesian optimization literature (Ozdemir *et al.*, 2025). EI quantifies the expected improvement over the best observed scalarized objective value by jointly considering the predictive mean and uncertainty of the surrogate model, thereby favouring designs that are both promising and informative. In contrast, UCB combines the predictive mean with a weighted uncertainty term, providing explicit control over the exploration–exploitation trade-off through a tunable parameter. The mathematical expressions for EI and UCB are given in Equations 18 and 19, respectively, and all symbols used are defined in the nomenclature.

$$\mathrm{EI}(x) = (s(x) - s^*)\Phi(z) + \sigma_s(x)\phi(z), \; z = \frac{s(x)-s^*}{\sigma_s(x)} \tag{18}$$

$$\mathrm{UCB}(x) = s(x) + \kappa\sigma_s(x) \tag{19}$$

As shown in Equation 19, the EI acquisition function consists of two contributions: a term proportional to the predicted improvement over the current best solution, which promotes exploitation, and a term proportional to the predictive uncertainty, which promotes exploration. Consequently, candidate points with high predicted performance are favoured when uncertainty is low, while points associated with higher uncertainty may be selected to improve model knowledge. In comparison, UCB offers a more direct and interpretable formulation of the exploration–exploitation trade-off through the parameter $\kappa$. Larger values of $\kappa$ encourage exploration of uncertain regions of the design space, whereas smaller values emphasize exploitation of regions with high predicted performance.

#### *2.3.3 Optimization of the acquisition function*

While the choice of acquisition function is critical, its practical application also requires an effective strategy for maximizing the acquisition function. In this study, acquisition function maximization was performed using a sampling-based optimization approach. A set of 2000 candidate mix designs was generated by uniformly sampling the design space while enforcing the physical and compositional constraints defined previously. Although sampling-based optimization does not guarantee identification of the global maximum of the acquisition function due to its reliance on a finite candidate set, it offers a simple and robust alternative to hybrid optimization methods. Völker *et al.* (2023) used a similar technique with a finite candidate pool to select feasible candidates from; albeit the candidate pool was not sampled but was considered based on an experimental database. In heavily constrained design spaces such as the one in this paper, sampling-based approaches are particularly well suited, as they naturally accommodate complex constraints and mixed variable types and often yield near-optimal solutions with modest computational effort. Increasing the number of points further could possibly result in a solution closer to the global optimum, but with increased computational memory demand. Hence, the decision to use 2000 points was based on a compromise between optimization performance and computation capacity.

For each candidate design, predictions from the surrogate models were obtained and the corresponding acquisition function value was evaluated. A post-sampling feasibility filtering step was then applied to enforce the performance, sustainability, and uncertainty-based constraints defined previously. Only candidates satisfying all post-sampling constraints were considered feasible.

Among the feasible candidates, the mix designs corresponding to the highest acquisition function values were selected and added to the dataset for subsequent model updates. Further details regarding the model update procedure are provided in Section 2.3.4.

#### *2.3.4 Database update*

With maximization of the acquisition function suggesting potential candidates to update the dataset, 3 different ways were considered for database update in this case which has been mentioned in this subsection.

- *Experimental database update with EI acquisition function*

In this case, two optimization iterations were conducted, in which mortar specimens were fabricated and tested in accordance with EN 196-1 (European Committee for Standardisation, 2016). At each iteration, the 2-day compressive strength was measured experimentally to enable rapid feedback within the optimization loop. In contrast, directly measuring 28-day compressive strength at each iteration would require a prolonged waiting period, rendering the sequential optimization process impractical. Therefore, the 28-day strength was incorporated through a model-based update, as defined in Equation 20, rather than direct experimentation. This strategy eliminates the need to wait 28 days for each iteration while still accounting for long-term performance. Furthermore, a safety factor was applied to the predicted 28-day strength by subtracting 1.96 times the model standard deviation (ensuring that the strength is higher than the lower bound of the 95 % confidence interval), thereby accounting for prediction uncertainty. This conservative formulation stabilizes decision-making in regions of high uncertainty and reflects a risk-averse design approach appropriate for material optimization under limited experimental data.

$$28 - day\ strength\ update = \mu_{28}(x) - 1.96\sigma_{28}(x) \tag{20}$$

- *Simulation-based database update with EI acquisition function*

In addition to experimental updates, a simulation-based database update was performed using the Expected Improvement (EI) acquisition function. Sequential iterations were considered until the Pareto dominance depth was zero which means the experimental loop seemed to converge to the global Pareto front. However, in this case, database updates were carried out using only the predicted mean responses ($\mu_2$ and $\mu_{28}$) obtained from the surrogate models. This simulation-based approach was adopted because conducting multiple sequential laboratory experiments was not feasible due to time and material constraints. The simulation-based updates therefore provide insight into the expected convergence behavior of the Bayesian optimization framework under idealized conditions.

- *Simulation-based database update with UCB acquisition function*

A similar simulation-based update strategy was employed using the Upper Confidence Bound (UCB) acquisition function. Iterations were performed (with a maximum of 10 iterations) with a starting $\kappa$ of 100 for the first few iterations and reduction to 20 for the next. This decreasing $\kappa$ schedule was adopted to emphasize exploration during the early stages of optimization, when model uncertainty is high, and to gradually shift toward exploitation as the surrogate models become more informative. In the present formulation, the predictive standard deviation of the 2-day compressive strength model, $\sigma_2(x)$, was normalized by the range of the 2-day strength values in the initial dataset during scalarization. Consequently, the chosen reduction of $\kappa$ from 100 to 20 corresponds to an effective reduction of approximately 5 to 1 in the exploration weight applied to the normalized uncertainty term. This scaling ensures that the relative influence of uncertainty in the UCB acquisition remains comparable to commonly reported $\kappa$ values in the BO literature (Ozdemir *et al.*, 2025).

Time or iteration-dependent adjustment of the exploration-exploitation parameter has been theoretically justified by Desautels, Krause and Burdick (2012), who proposed adaptive $\kappa$ schedules that account for the evolution of predictive uncertainty as the model is updated. In applied settings, Ozdemir *et al.* (2025) employed $\kappa$ values in the range of 0.5–3 to toggle between exploratory and exploitative regimes using UCB. As with the EI-based simulations, database updates in the UCB case were performed using predicted mean responses rather than experimentally measured values.

#### *2.3.5 Baseline comparison*

The results obtained from both experimental and simulation-based Bayesian optimization updates using EI and UCB were benchmarked against the Pareto front constructed from the full 51-mix dataset. For consistency, only mix designs with a water-to-binder ratio of at least 0.48 were considered when constructing the reference Pareto front, ensuring fair comparison across similar flowability conditions. This baseline comparison aims to assess the ability of the BO framework to converge toward the global Pareto front using a substantially reduced number of evaluations, starting from a carefully selected initial dataset comprising approximately 20% of the full dataset.

## 3. Results and Discussions

### *3.1 Kernel performance and selection*

Based on the grid search-based length-scale optimization of kernels, the radial basis function (RBF) kernel provided the best overall predictive performance among the kernels considered. The length-scale combinations suggested by the grid search, along with the corresponding train and test $R^2$values, are summarized in Table 4.

**Table 4**: Grid-search length-scale optimization results for 2-day and 28-day compressive strength GP models

| Length scale | Length scale combination* | Train $R^2$ | Test $R^2$ |
|---|---|---|---|
| 2-day compressive strength GP model | (20, 20, 0.1, 3, 3) | 0.97 | 0.75 |
| 28-day compressive strength GP model | (20, 20, 0.5, 2, 2) | 0.95 | 0.38 |

*Length scale combination = (copper slag, limestone, water-to-binder ratio, cement type, Fe-rich slag fineness)

For the 2-day compressive strength GP model, the grid search identified a length-scale combination of (20, 20, 0.1, 3, 3), yielding a training $R^2$of 0.97 and a testing $R^2$of 0.75. However, a slightly modified length-scale combination of (20, 20, 0.05, 3, 3) was ultimately selected, as reducing the length scale associated with the water-to-binder ratio improved model sensitivity while maintaining comparable generalization performance (train $R^2 = 0.98$, test $R^2 = 0.74$). This adjustment was made to avoid overly smooth response surfaces that could limit the effectiveness of uncertainty-driven exploration in subsequent BO iterations.

For the 28-day compressive strength GP model, the selected RBF kernel (Table ) achieved a training $R^2$of 0.95 but a lower testing $R^2$of 0.38, reflecting the increased variability and the possibility of Fe-rich slag showing late-age reactivity influencing model behaviour at 28 days (Hallet *et al.*, 2023). The 28-day compressive strength GP model was used mainly as a constraint in the BO loop in this paper.

In contrast to the non-linear kernels, the linear kernel showed better performance in terms of test $R^2$ with values of 0.80 and 0.60 for 2-day and 28-day compressive strength GP models respectively. The train $R^2$ values were 0.90 and 0.86 for 2-day and 28-day compressive strength GP models respectively. While the linear kernel achieved superior test $R^2$ (0.80 vs 0.75), RBF was selected for the BO loop because it provides more conservative uncertainty estimates in regions far from training data. In active learning, where the goal is sequential exploration, the quality of uncertainty quantification is as important as mean prediction accuracy. Linear kernels produce uncertainty estimates that reflect distance from the origin rather than proximity to observed data, making them insensitive to the local data structure. In contrast, the RBF kernel appropriately increases posterior uncertainty in regions sparse of observations, leading to a more principled exploration-exploitation balance in the acquisition function.

Overall, kernel selection prioritized stable generalization behaviour and physically meaningful sensitivity over maximum training accuracy, consistent with the objectives of data-efficient experimental design in small-sample cementitious systems. The 28-day compressive strength model shows greater signs of overfitting than the 2-day compressive strength model which can be seen from the higher difference in the train and test $R^2$ values in the case of the 28-day compressive strength model as compared to the 2-day compressive strength model.

### *3.2 Experimental database update with EI*

#### *3.2.1 First iteration*

Using $D_0$ as the initial database, the first iteration of the BO loop was performed using the Expected Improvement (EI) acquisition function. Table 5 lists the candidate mix designs with the highest EI values for each combination of Fe-rich slag fineness (2700 and 4400 $cm^2/g$) and cement type (CEM I 52.5N and R). Figure 1 shows the positions of the four selected points (as per EI per combination) along with all the uniformly sampled points and the feasible points satisfying the criteria mentioned in section 2.2.3. Among these candidates, the mix corresponding to the 2700-N combination was selected for laboratory testing, as the alternative candidate requiring slag with a fineness of 4400 $cm^2/g$ could not be prepared due to material availability constraints.

**Table 5**: Predicted mix designs from first iteration of the BO loop

| Categorical variable combination | Fe-rich slag (wt%) | Limestone (wt%) | Portland cement (wt%) | Water-to-binder ratio |
|---|---|---|---|---|
| 4400-N | 24.20 | 24.50 | 51.30 | 0.48 |
| 2700-N | 21.90 | 21.90 | 56.20 | 0.48 |
| 2700-R | 26.40 | 20.00 | 53.60 | 0.48 |
| 4400-R | 28.50 | 19.40 | 52.00 | 0.48 |

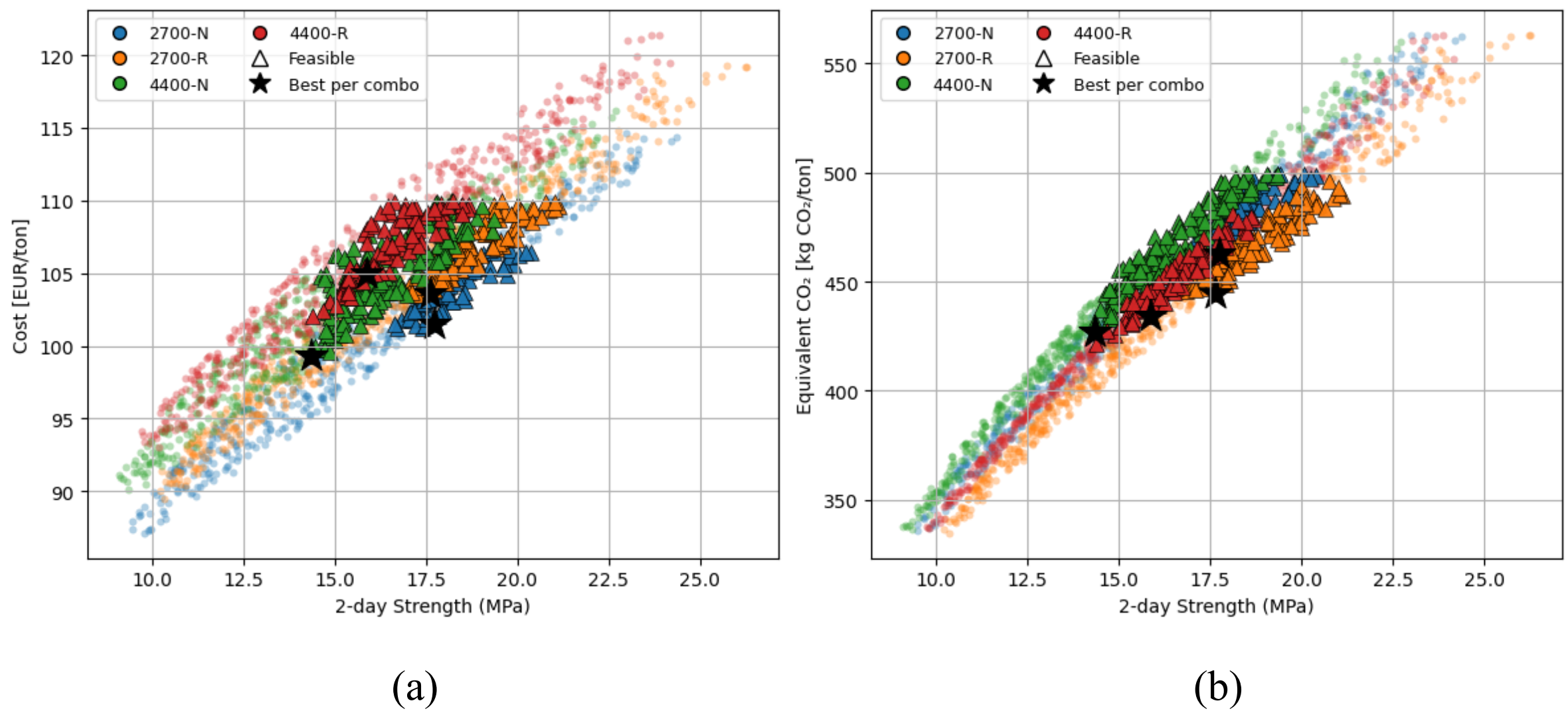

**Figure 1**: (a) Cost vs 2-day strength uniformly sampled points (satisfying box constraints); (b) Equivalent $CO_2$ vs 2-day strength uniformly sampled points (satisfying box constraints)

Figure 1 illustrates the uniformly sampled candidate designs in the two-dimensional objective projections of cost vs 2-day compressive strength and equivalent $CO_2$ emissions vs 2-day compressive strength. The plots clearly show that, while the design space is broadly sampled, only a relatively narrow region satisfies the imposed feasibility constraints. This concentration of feasible points highlights the highly constrained nature of the optimization problem, where constraints on compressive strength, cost, $CO_2$ emissions, and mix proportions significantly reduce the effective search space available to the acquisition function. As a result, the probability of identifying improved candidate designs, regardless of the specific acquisition function employed, is inherently limited, since most sampled points are rendered infeasible prior to evaluation. In contrast, in the absence of such constraints, a substantially larger fraction of the design space would be available for exploration, thereby increasing the likelihood of identifying improved solutions. The effect becomes even more pronounced for stricter performance requirements, such as higher cement strength classes (e.g., 42.5N or 52.5N), where the feasible region is expected to shrink further. Gardner, Jake and Kusner (2014) report such issues for constrained BO problems and recommend the use of specified constraint-aware acquisition functions for such problems. In such cases, uniform sampling might not be the best solution. However, in this paper, the uniform sampling was conducted only on points satisfying the material constraints, increasing the probability of finding points that satisfy the strength, cost and equivalent $CO_2$ inequality constraints as well.

Table 6 reports the Gaussian process (GP) predictions for both 2-day and 28-day compressive strength, including the predictive mean and standard deviation, alongside the experimentally measured values.

**Table 6**: 2-day and 28-day compressive strength GP model predictions and experimental values

| Categorical variable combination | 2-day compressive strength GP predictions $\mu_2(x) \pm \sigma_2(x)$ | 28-day compressive strength GP predictions $\mu_{28}(x) \pm \sigma_{28}(x)$ | 2-day compressive strength experimental results | 28-day compressive strength experimental results |
|---|---|---|---|---|
| 2700-N | 17.74 ± 2.28 | 40.05 ± 3.52 | 13.70 ± 0.28 | 33.00 ± 0.60 |

As expected in the first iteration with a limited initial dataset, the GP predictions exhibit a degree of overconfidence relative to the experimental results. Specifically, the experimentally measured 28-day compressive strength falls slightly below the lower bound of the 95% confidence interval, while for the 2-day compressive strength the experimental value lies within the confidence interval, with a deviation of approximately 0.6% from its lower bound.

Such discrepancies are common in early iterations of Bayesian optimization, where surrogate models are trained on sparse data and local variability has not yet been fully learned. Importantly, the Bayesian optimization framework explicitly accounts for this uncertainty through conservative update rules. For the subsequent iteration, the experimentally measured 2-day compressive strength was used to update the corresponding GP model, while the lower bound of the 95% confidence interval from model prediction ($\mu_{28}(x) - 1.96\sigma_{28}(x)$) was employed for updating the 28-day compressive strength GP model, as described in Section 2.2.5. This update strategy mitigates the impact of model overconfidence and promotes robust convergence of the optimization process in later iterations. Such a strategy was implemented by Pfeiffer *et al.* (2024) where higher early-age concrete mixes were predicted by GP models with uncertainty incorporation. Additionally, the cost and $CO_2$ of the reported mix design were 101.5 EUR/ton and 462.5 kg $CO_2$/ton respectively. The equivalent $CO_2$ value is lower than the 2050 target set by CEMBUREAU (CEMBUREAU, 2020), albeit with two strength classes lower.

#### *3.2.2 Second iteration*

Using results from the first iteration (with an updated database $D_1$ with 11 data points), the second iteration of the BO loop was conducted. Like section 3.2.1, Table 7 lists the candidate mix designs with the highest EI values for each combination of cement type and Fe-rich slag fineness. Among these candidates, the mix corresponding to the 2700-R combination was selected for laboratory testing, as the alternative candidate requiring slag with a fineness of 4400 cm²/g could not be prepared due to material availability constraints.

**Table 7**: Predicted mix designs from second iteration of the BO loop

| Categorical variable combination | Fe-rich slag (wt%) | Limestone (wt%) | Portland cement (wt%) | Water-to-binder ratio |
|---|---|---|---|---|
| 4400-R | 32.20 | 16.10 | 51.70 | 0.48 |

| 4400-N | 29.80 | 16.10 | 54.10 | 0.50 |
|---|---|---|---|---|
| 2700-R | 29.30 | 15.30 | 55.40 | 0.48 |
| 2700-N | 29.20 | 12.30 | 58.50 | 0.48 |

Table 8 reports the Gaussian process (GP) predictions for both 2-day and 28-day compressive strength at the second iteration of the Bayesian optimization loop, including the predictive mean and associated standard deviation, alongside the experimentally measured values. Compared to the first iteration, an improvement in predictive performance is observed for both early-age and later-age strength models, in terms of prediction uncertainty.

**Table 8**: 2-day and 28-day compressive strength GP model predictions and experimental values

| **Categorical variable combination** | **2-day compressive strength GP predictions** $\mu_2(x) \pm \sigma_2(x)$ | **28-day compressive strength GP predictions** $\mu_{28}(x) \pm \sigma_{28}(x)$ | **2-day compressive strength experimental results** | **28-day compressive strength experimental results** |
|---|---|---|---|---|
| 2700-N | 17.64 ± 1.65 | 37.80 ± 2.37 | 16.06 ± 0.43 | 36.12 ± 2.78 |

In particular, the experimentally measured 2-day and 28-day compressive strength values lie above the lower bounds of the corresponding 95% confidence intervals predicted by the GP models. This indicates that the surrogate models provide more conservative and reliable uncertainty estimates following the database update performed after the first iteration. The improved agreement between predictions and experimental results reflects the effect of incorporating new information into the GP models, thereby reducing model overconfidence and enhancing predictive robustness. In addition to meeting the compressive strength requirements, the selected mix design also satisfies the imposed sustainability constraints, with a cost of 105.4 EUR/ton and an equivalent $CO_2$ emission of 483.5 kg $CO_2$/ton. Similar to results from the first iteration, the equivalent $CO_2$ of the predicted blended cement mix satisfies the 2050 target proposed by CEMBUREAU (CEMBUREAU, 2020). These results demonstrate that the BO framework can progressively identify mix designs that simultaneously satisfy performance, economic, and environmental criteria as the sequential learning process advances.

Furthermore, the points generated for the two iterations were compared to the global Pareto front with the entire dataset (with only water-to-binder ratio ≥ 0.48). Figure 2 show the two faces of the 3-objective (2-day compressive strength, cost and equivalent $CO_2$) Pareto front of the points and shows the position of the points from the two iterations (in red squares).

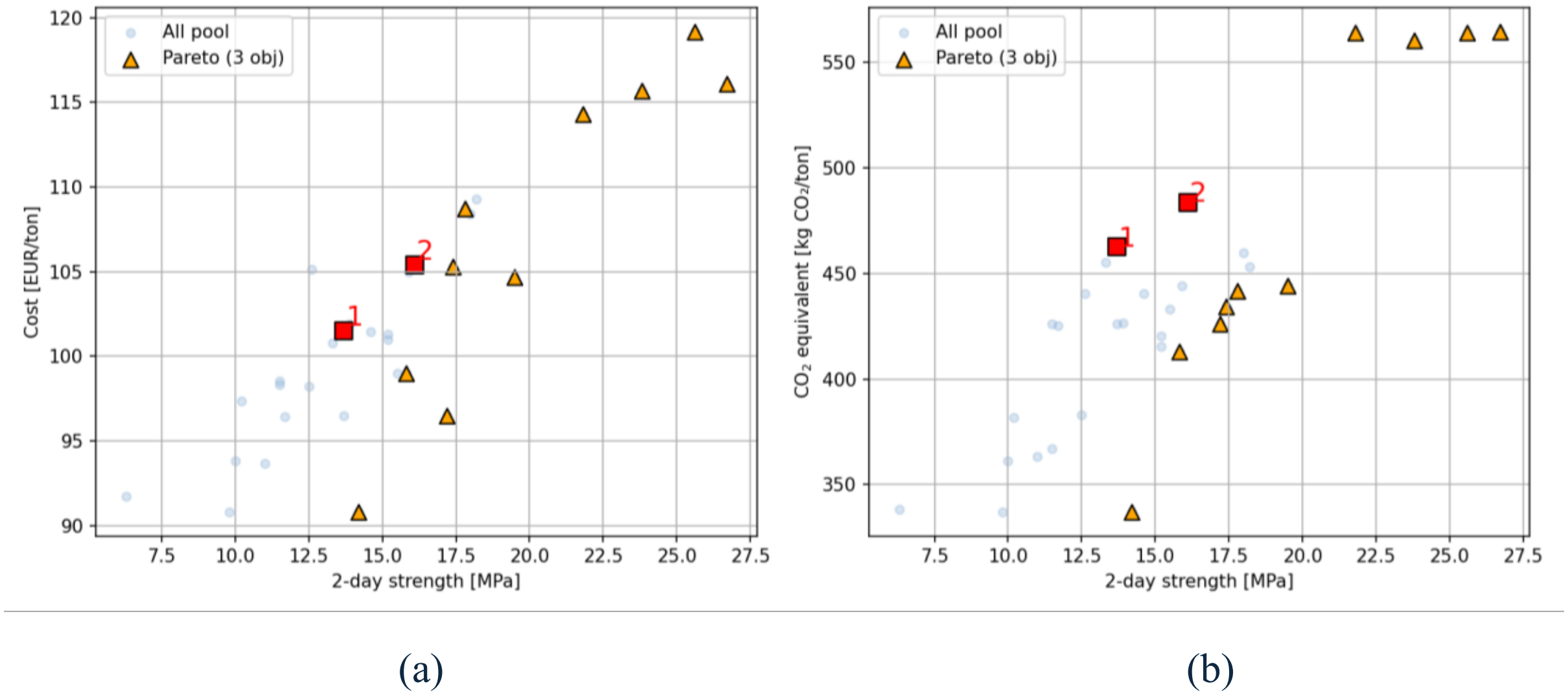

(a) (b)

**Figure 2**: Two faces (a): Cost vs 2-day strength and (b): $CO_2$ equivalent vs 2-day strength of the 3-objective global Pareto front with points generated by the 2 iterations and two faces

Figure 2 (a and b) illustrates the position of the mix designs generated in the first two BO iterations relative to the global three-objective Pareto front constructed from the full dataset. In both objective projections, the candidate selected in the second iteration is visibly closer to the non-dominated region than the candidate obtained in the first iteration, indicating progressive movement toward the Pareto front. This trend is further quantified through dominance analysis with respect to the global Pareto front. The mix design generated in the first iteration is dominated by eight Pareto-optimal points, whereas the design generated in the second iteration is dominated by only three points. The reduction in dominance depth demonstrates that the BO framework is effectively guiding the search toward improved trade-offs among compressive strength, cost, and $CO_2$ emissions. While full convergence to the global Pareto front cannot be guaranteed within a limited number of iterations, along with a conservative model update strategy adopted in this study, the observed reduction in dominance indicates meaningful directional convergence. Notably, a solution dominated by only three Pareto-optimal points was identified using an initial dataset of ten points and a total of twelve evaluations, whereas the reference Pareto front was constructed from 51 mix designs (restricted to water-to-binder ratios $\geq 0.48$, so 30 points here). This comparison highlights the ability of BO to achieve substantial progress towards Pareto-optimal solutions with significantly fewer evaluations through sequential, data-efficient learning.

Several data points in the pool achieve higher 2-day compressive strength at lower cost or $CO_2$ equivalent but do so at wb ratios of 0.45, which lie outside the feasible design space imposed in the BO loop (wb $\geq$ 0.48). As a result, direct comparisons between BO-generated candidates and training points with wb = 0.45 are not strictly fair, since the latter violate the imposed workability constraint. Consequently, the BO algorithm was effectively optimizing within a reduced and more practically relevant design space, where the apparent performance gap reflects constraint enforcement rather than a failure of the optimization strategy. This observation highlights that, particularly for small datasets, the consistency and representativeness of the training data with respect to the imposed feasibility constraints play a critical role in the effectiveness of BO-based optimization. In such cases, data quality and constraint alignment can be more influential than the choice of acquisition function in determining observed optimization performance.

*3.3 Simulation-based database update with EI*

This section presents results of iterations with simulated data (with $\mu_2(x)$ and $\mu_{28}(x)$ as database updates) to demonstrate the performance of the GP-based BO loop with EI as an acquisition function with a higher number of iterations as compared to the experimental update in section 3.2. Table 9 shows the predicted mix design with the highest EI at each iteration while Table 10 shows the corresponding 2-day compressive strength, 28-day compressive strength, cost, $CO_2$ equivalent and Pareto dominance depth of the different iterations.

**Table 9**: Predicted mix design with the highest EI at each iteration

| Iteration no. | Fe-rich slag (wt%) | Limestone (wt%) | Cement (wt%) | Water-to-binder ratio | Cement type | Fe-rich slag fineness |
|---|---|---|---|---|---|---|
| 1 | 24.20 | 24.50 | 51.30 | 0.48 | N | 4400 |
| 2 | 25.50 | 26.10 | 48.40 | 0.48 | N | 4400 |
| 3 | 21.90 | 23.20 | 54.90 | 0.50 | N | 2700 |
| 4 | 20.40 | 26.60 | 53.00 | 0.48 | N | 2700 |
| 5 | 22.70 | 27.10 | 50.30 | 0.50 | N | 2700 |
| 6 | 23.80 | 26.70 | 49.50 | 0.48 | N | 2700 |

**Table 10**: 2-day compressive strength, 28-day compressive strength, cost, CO2 equivalent and Pareto dominance depth of the different iterations

| Iteration no. | 2-day compressive strength [MPa] | 28-day compressive strength [MPa] | Cost [EUR/ton] | $CO_2$ equivalent [kg $CO_2$/ton] | Pareto dominance depth (2-day strength, cost and $CO_2$ equivalent) |
|---|---|---|---|---|---|
| 1 | 14.35 ± 2.21 | 38.87 ± 2.92 | 99.3 | 426.7 | 2 |
| 2 | 13.07 ± 1.49 | 37.13 ± 2.17 | 96.8 | 405.1 | 1 |
| 3 | 16.10 ± 1.97 | 39.07 ± 3.29 | 100.2 | 452.4 | 1 |
| 4 | 15.85 ± 1.66 | 38.18 ± 2.55 | 98.0 | 436.9 | 1 |
| 5 | 14.17 ± 1.21 | 36.22 ± 1.84 | 95.7 | 417.1 | 1 |
| 6 | 14.31 ± 1.03 | 35.69 ± 1.47 | 95.2 | 411.5 | 0 |

From Table 9, six BO iterations were simulated using database updates based on the predictive means $\mu_2(x)$ and $\mu_{28}(x)$. A practical stopping criterion was defined as the attainment of a Pareto dominance depth of zero, indicating that the selected mix design is non-dominated with respect to the global Pareto front constructed from the full dataset. As shown in Table 10, this criterion was satisfied at iteration 6. It is important to note that achieving a Pareto dominance depth of zero does not guarantee convergence to the true global Pareto front. Additional iterations could, in principle, further refine trade-offs or explore nearby regions of the Pareto front, as reported in previous BO studies (Ozdemir *et al.*, 2025). However, in heavily constrained design spaces such as the one considered in this study, the likelihood of substantial further improvement diminishes rapidly, and conducting additional experiments may be impractical due to cost and time constraints.

Furthermore, a gradual reduction in predictive uncertainty, reflected by decreasing $\sigma_2(x)$ and $\sigma_{28}(x)$, is observed across iterations. This behaviour is consistent with established Bayesian optimization theory, wherein Gaussian process posterior uncertainty contracts locally as new data are incorporated (Hernández-Lobato *et al.*, 2015). This indicates that the surrogate models dynamically learn from the updated dataset at each iteration.

A comparison between different update strategies further highlights the importance of uncertainty-aware model updates. The Pareto dominance depths for the first two iterations were eight and three, respectively, in the case of experimental and conservative model updates (section 3.2), whereas they were two and one in the case of mean-based ($\mu_2(x)$ and $\mu_{28}(x)$) updates (section 3.3). This suggests that GP predictions based solely on mean values can be overconfident, leading to solutions that appear closer to the Pareto front than those obtained through uncertainty-aware updates. Consequently, incorporating experimental data and uncertainty considerations is expected to yield designs that more reliably approach the global Pareto front. Figure 3 shows the position of the six points (unordered) compared to the global Pareto front with three objectives (2-day compressive strength, cost and equivalent $CO_2$) but shown in two faces.

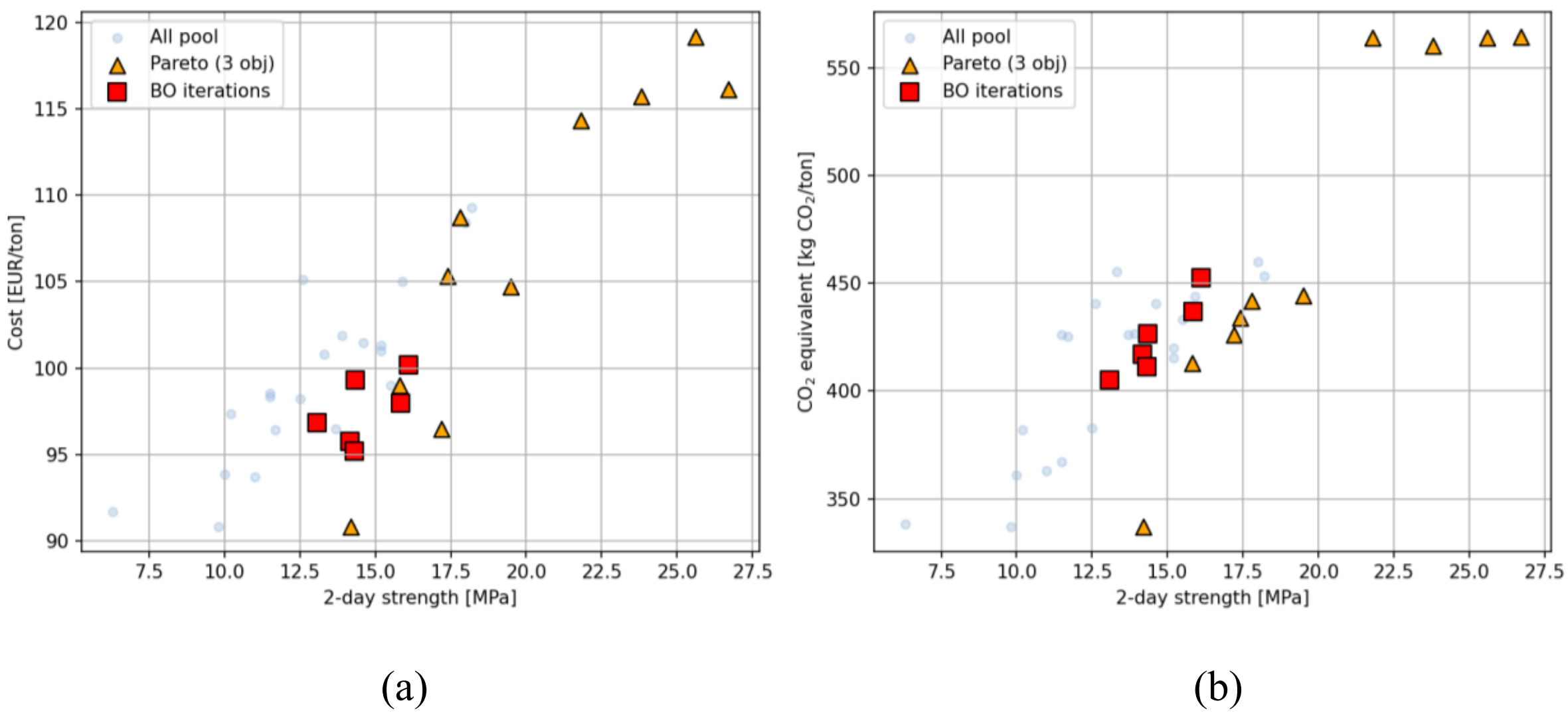

**Figure 3**: (a) Cost vs 2-day strength with simulated EI-based BO loop; (b) $CO_2$ equivalent vs 2-day strength with simulated EI-based BO loop

From Figure 3, it can be observed that the simulated EI-based BO iterations cluster closely in objective space. This clustering behaviour is particularly pronounced due to the heavily constrained nature of the design space considered in this study, which restricts feasible solutions to narrow regions of the cost-strength-$CO_2$ equivalent trade-off surface. Additionally, the EI acquisition function is known to place relatively greater emphasis on exploitation (Shoyeb Raihan *et al.*, 2024), especially as predictive uncertainty decreases. As a result, EI tends to select candidate designs in regions that are already predicted to perform well, further contributing to the observed clustering of successive iterations.

It is also important to note that, in the present formulation, predictive uncertainties governing the exploratory behaviour of the acquisition function arise solely from the Gaussian process models for 2-day and 28-day compressive strength ($\sigma_2(x)$ and $\sigma_{28}(x)$). In contrast, cost and $CO_2$ equivalent are modelled deterministically and therefore do not contribute to uncertainty in the acquisition function. These further bias the optimization process toward exploitation-dominated behaviour, leading to closely spaced BO iterations in the objective space.

### *3.4 Simulation-based update with UCB*

This section presents results of 10 iterations (with first five iterations with $\kappa$ as 100 and next five with $\kappa$ as 20) with simulated data (with $\mu_2(x)$ and $\mu_{28}(x)$ as database updates) to demonstrate the performance of the GP-based BO loop with UCB as an acquisition function with a higher number of iterations. Compared to section 3.3, the stopping criterion was strictly defined at 10 iterations simulating real life experimental budget constraints. Table 11 shows the predicted optimal mix design at each iteration with their corresponding $\kappa$ value while Table 12 shows the predicted 2-day and 28-day compressive strength values from their GP models (along with uncertainties) and their cost, $CO_2$ equivalent and Pareto dominance depth corresponding to the original dataset (with water-to-binder ≥ 0.48).

**Table 11**: Predicted optimal mix designs at each iteration based on UCB acquisition function

| Iteration | $\kappa$ | Fe-rich slag (wt%) | Limestone (wt%) | Cement (wt%) | Water-to-binder ratio | Cement type | Fe-rich slag fineness ($cm^2/g$) |
|---|---|---|---|---|---|---|---|
| 1 | 100 | 36.8 | 9.1 | 54.1 | 0.49 | R | 4400 |
| 2 | 100 | 20.3 | 23.3 | 56.4 | 0.50 | N | 2700 |
| 3 | 100 | 29.2 | 19.3 | 51.6 | 0.48 | N | 4400 |
| 4 | 100 | 21.4 | 24.0 | 54.6 | 0.48 | R | 2700 |
| 5 | 100 | 20.9 | 27.7 | 51.5 | 0.48 | N | 2700 |
| 6 | 20 | 25.0 | 26.0 | 49.0 | 0.50 | N | 4400 |
| 7 | 20 | 21.0 | 27.8 | 51.2 | 0.49 | R | 2700 |
| 8 | 20 | 32.7 | 19.2 | 48.1 | 0.50 | N | 4400 |
| 9 | 20 | 24.9 | 25.3 | 49.8 | 0.50 | R | 2700 |
| 10 | 20 | 30.1 | 22.3 | 47.6 | 0.48 | R | 4400 |

**Table 12**: Predicted 2-day and 28-day compressive strength (with uncertainties), cost and $CO_2$ equivalent and corresponding Pareto dominance depth

| Iteration no. | 2-day compressive strength [MPa] | 28-day compressive strength [MPa] | Cost [EUR/ton] | $CO_2$ equivalent [kg $CO_2$/ton] | Pareto dominance depth (2-day strength, cost and $CO_2$ equivalent) |
|---|---|---|---|---|---|
| 1 | 16.09 ± 2.64 | 41.23 ± 4.32 | 109.8 | 454.8 | 5 |
| 2 | 16.80 ± 2.23 | 40.16 ± 3.59 | 101.3 | 463.0 | 1 |
| 3 | 14.55 ± 2.05 | 39.52 ± 3.15 | 101.1 | 432.0 | 3 |
| 4 | 17.88 ± 1.70 | 39.56 ± 2.64 | 103.5 | 449.4 | 0 |
| 5 | 15.62 ± 1.61 | 37.73 ± 2.42 | 96.6 | 425.6 | 0 |
| 6 | 13.09 ± 1.48 | 37.29 ± 1.87 | 97.3 | 409.9 | 1 |
| 7 | 15.71 ± 1.37 | 37.54 ± 2.42 | 99.9 | 423.3 | 1 |
| 8 | 12.54 ± 1.28 | 37.03 ± 1.75 | 98.8 | 406.5 | 1 |
| 9 | 14.86 ± 1.13 | 36.04 ± 1.71 | 99.2 | 414.4 | 1 |
| 10 | 13.40 ± 1.20 | 36.17 ± 1.83 | 100.8 | 401.2 | 1 |

From Table 11, it can be observed that the variation of Fe-rich slag, limestone, and Portland cement contents across iterations is non-linear and does not follow a monotonic trend. This behaviour is expected in BO, as the selected mix designs result from a balance between predicted performance and uncertainty rather than from a deterministic optimization path. The observed trends reflect the combined influence of raw material proportions on 2-day and 28-day compressive strength, cost, and $CO_2$ equivalent. It is further observed that mix designs with Fe-rich slag contents exceeding approximately 30 wt% are frequently Pareto dominated within the considered multi-objective setting. This suggests that, under the imposed constraints and within the explored design space, the use of higher Fe-rich slag contents offers limited trade-off benefits. A general reduction in predictive uncertainty, reflected by decreases in $\sigma_2(x)$ and $\sigma_{28}(x)$, is observed from the first to the tenth iteration, although not strictly monotonically. Such behaviour is characteristic of BO, where surrogate models progressively improve as additional data are incorporated.

The key distinction between the UCB-based updates presented here and the EI-based updates discussed in section 3.3 lies in the distribution of Pareto dominance depth across iterations. When $\kappa$ was set to 100 (strongly exploratory regime) during the first five iterations, the dominance depths spanned a wide range (from zero to five), indicating active exploration of diverse regions of the constrained design space. In contrast, after $\kappa$ was reduced to 20 (more exploitative regime), the dominance depth of selected points decreased and stabilized to one. In this case, the value of κ was heuristically changed to suit the problem; however, there are theoretical techniques in literature that govern the selection of $\kappa$ in the problem (Srinivas *et al.*, 2010). Although strict convergence of the BO loop cannot be guaranteed with a limited experimental budget, the results demonstrate that the UCB-based strategy was able to identify Pareto non-dominated solutions while maintaining exploratory behaviour during the early iterations.

Figure 4 shows the distribution of the points sampled by BO iterations along with the Pareto front of the master data set.

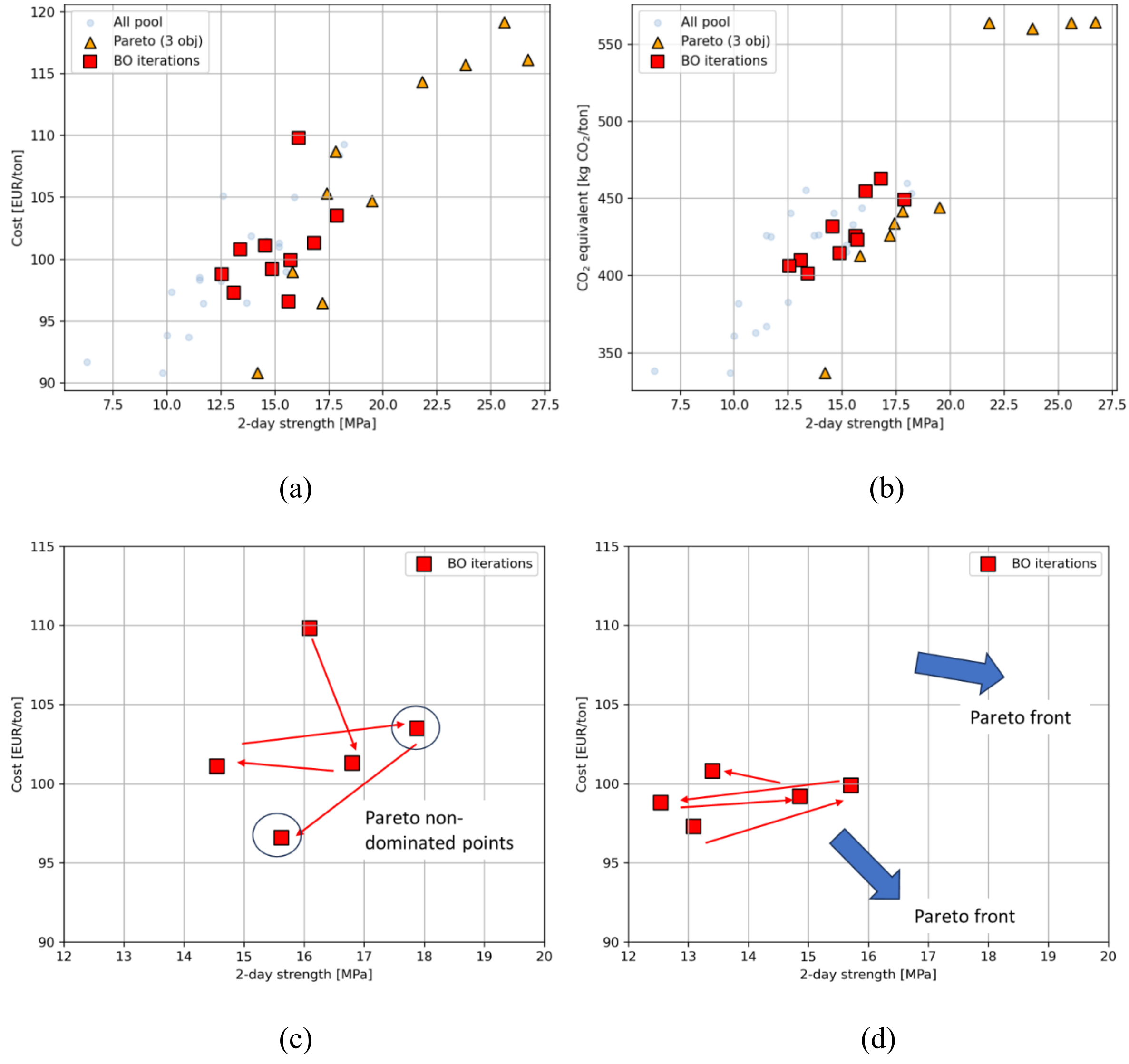

**Figure 4**: UCB-based BO sampled points (simulated): (a) Cost vs 2-day strength; (b) $CO_2$ equivalent vs 2-day strength; (c) Cost vs 2-day strength for exploration-based sampling (first five points); (d) Cost vs 2-day strength for exploitation-based sampling (last five points)

From Figure 4, the BO-sampled points obtained using the UCB acquisition function appear more scattered compared to those sampled using EI-based BO (Figure 3). This indicates a stronger exploration-exploitation trade-off inherent to UCB-based sampling, where exploration is explicitly encouraged through the uncertainty term in the acquisition function. In contrast, EI-based sampling tends to focus more on exploitation by prioritizing regions with higher expected improvements. The increased dispersion of the UCB-sampled points suggests that the algorithm actively explores a broader region of the design space, potentially improving global search capabilities. From **Figure *4*** (c) and (d), it can be observed that the Bayesian optimization trajectory shows a general tendency toward Pareto non-dominance; however, this progression is not monotonic, particularly during the exploration regime (**Figure *4***c). For clarity, only the cost versus 2-day compressive strength projections is shown in **Figure *4*** (c) and (d). The corresponding $CO_2$ equivalent versus 2-day strength projections exhibit significant point clustering, primarily because CEM I 52.5N and CEM I 52.5R were assumed to have identical $CO_2$ equivalents, limiting visual separability in that objective space.

However, the relaxation of the constraints may have further amplified this exploratory behaviour, allowing the UCB-based BO to sample points in more promising regions of the objective space. As a result, a more diverse set of candidate solutions could have been generated, potentially leading to more promising Pareto-optimal points.

## 4. Conclusions

This study investigated the applicability of BO as a data-efficient framework for multi-objective blended cement mix design under practical experimental and compositional constraints. The conclusions are structured to directly address the research questions posed in this work.

- *Effectiveness of BO under limited data*

BO successfully identified Pareto non-dominated mix designs using a very small initial dataset of ten experiments, demonstrating its suitability for laboratory-scale cement optimization where experimental budgets are severely constrained. Contrary to the common perception that AI-assisted methods require thousands of data points which is an unrealistic requirement for cement and concrete research, this study shows that meaningful optimization can be achieved with a data-to-independent-variable ratio as low as 2:1, provided that the BO framework is carefully designed. In particular, the informed selection of Gaussian process (GP) models, kernels, feature-specific length scales, and acquisition functions proved critical. Additionally, all data were generated within the same laboratory by the same operator, thereby minimizing user-dependent variability that commonly affects studies relying on heterogeneous literature datasets and can otherwise degrade BO performance.

In the context of categorical variables, they were one-hot encoded in this paper and then used in the problem. No special consideration is necessary for such variables in this case, but it is indeed possible to use special kernels.

- *Performance of GP surrogate models*

The predictive performance of the GP surrogate models differed markedly between early-age and late-age compressive strength. For 2-day compressive strength, test $R^2$ values up to 0.75 were achieved, whereas the corresponding maximum test $R^2$ for the 28-day compressive strength model was limited to 0.38, despite both models exhibiting training $R^2$ values exceeding 0.90. While this discrepancy can partly be attributed to the small training dataset, it also suggests a higher intrinsic variability in the 28-day strength data, potentially arising from late-age reactions of Fe-rich slag. Furthermore, kernel selection and feature-dependent length scales were determined using a brute-force grid search, which, although transparent, may not yield globally optimal hyperparameters and could contribute to sub-optimal generalization performance. Additionally, the RBF kernels selected in this paper failed to perform better than their linear counterpart in terms of predictive ability. However, the RBF kernel was still considered owing to the chance of improved performance with iterations and the ability to use different length scales for different independent variables. Also, RBF kernels perform better than linear kernels in the case of uncertainty quantification which this paper considers a key element.

- *Acquisition functions and database updates*

Expected Improvement (EI) and Upper Confidence Bound (UCB) were employed as acquisition functions, and their comparative performance in the context of blended cement optimization constitutes a novel contribution of this work. In experimental BO updates using EI (limited to two iterations due to time and budget constraints), 2-day compressive strength was updated using measured experimental values, while 28-day strength was conservatively updated using the lower bound of the 95% confidence interval. This risk-averse strategy effectively bypasses the need to wait 28 days per iteration, addressing a critical bottleneck in applying BO to cement and concrete systems. However, this could possibly result in underestimated predictions of mix designs especially related to 28-day compressive strength. This is a cost the user pays for quicker mix design optimization and not waiting 28 days for an iteration. In addition, simulation-based BO loops were implemented using EI and UCB, with the latter employing different exploration parameters ($\kappa = 100$ and 20), showing the differences between cases where exploration and exploitation are preferred, respectively.

- *Cement performance targets and application-oriented design*

Starting from an initial dataset of ten experiments and applying BO-guided updates, blended cement mix designs meeting the 32.5N strength class were successfully generated, with two designs experimentally validated. All identified solutions satisfied a cost constraint of 110 EUR/ton and exhibited $CO_2$ equivalent emissions below 500 kg $CO_2$/ton, consistent with the long-term decarbonization targets proposed by CEMBUREAU for cementitious materials, albeit for a lower strength class.

- *Case-wise decision-making implications*

Finally, the comparative analysis of EI- and UCB-based BO strategies highlights that no single acquisition function is universally optimal. Instead, the choice of BO strategy should be guided by the specific time, budget, and risk constraints faced by the experimenter. EI is well-suited for exploitative, cost-efficient refinement near promising regions, whereas UCB offers greater flexibility for exploration under high uncertainty. These insights provide practical guidance for deploying BO as a decision-support tool in cement and concrete research and is a clear improvement to the previous works of the authors where a traditional experiment-analyse-predict approach was considered (Som *et al.*, 2026).

## 5. Recommended directions for future work

There are several areas that needs further work in the context of applying BO for cement and concrete research. They have been mentioned in the points below which are classified broadly into four categories.

1. *GP models and kernels*: GP models tend to work well as a surrogate model in BO and has been used extensively in BO literature. However, the kernels used are extremely critical in ensuring a stable performing GP model useful as a surrogate in the BO loop. In this context, categorical features should be considered differently as opposed to the current and previous works. Special kernels like Hamming kernels (Qiu *et al.*, 2024) can possibly be used for categorical features. Furthermore, non-stationary combined kernels have also shown promise in other fields (Ozdemir *et al.*, 2025), the knowledge of which can be transferred to the field of cement and concrete mix design optimization.

2. *Acquisition function optimization*: In place of sampling points and selecting on the basis of maximized acquisition functions, the acquisition function can itself be maximized using a hybrid solver (Noack and Funke, 2017) which combines a local and a global solver which can efficiently negotiate multiple local minima without taking a long time to converge. This is expected to increase the chances of finding a global optimum in the BO loop.
3. *Practical applications*: The techniques (or similar techniques) discussed in this paper should be applied more in research in place of traditional techniques, wherever possible. This will eventually help in benchmarking techniques and identifying best practices in the context of a powerful technique like BO for cement and concrete research. There are automated BO and experimental design tools used like BayBE (Fitzner *et al.*, 2025) which can perform tasks similar to the one shown in this paper. However, care needs to be taken when applying such tools, especially when employing black box models.
4. *Objective functions:* In this study, the focus was primarily on compressive strength, cost and $CO_2$ equivalent. However, important durability characteristics like freeze-thaw, acid attack and others should be included in the process of decision-making. However, in those cases, testing times are usually long (up to 12 months at times) which defeats the purpose of BO which is primarily used for accelerated material discovery. In that case, possible rapid alternatives or high-throughput techniques can be a solution to include such properties in a BO loop.

**Funding**

This project has received funding from the European Union's Horizon 2020 research and innovation programme under grant agreement No 958307. The project website is https://h2020harare.eu/.

**Declarations**

*Data availability:* The data used is available in the cited Zenodo link. The codes used in this paper will be publicly released upon acceptance of this paper.

*Conflicts of interest*: The authors state that there is no conflict of interest.

*Use of GenAI*: During the preparation of this work the author(s) used ChatGPT and Claude to polish the already written text for better textual flow and to generate the graphical abstract. After using this tool/service, the author(s) reviewed and edited the content as needed and take(s) full responsibility for the content of the published article. ChatGPT was also used for occasional debugging of code snippets.